\documentclass[a4paper,12pt,oneside]{article}
\usepackage{amsmath,amsfonts,amssymb,amsmath,latexsym,amsthm,enumerate}

\newtheorem{thm}{Theorem}

\theoremstyle{definition}

\theoremstyle{plain}

\begin{document}
\title {Some identities arising from Sheffer sequences for the powers of Sheffer pairs under umbral composition.}
\author{by \\Dae San Kim and Taekyun Kim}\date{}\maketitle

\begin{abstract}
\noindent In this paper, we study some properties of Sheffer sequences for the powers of Sheffer pairs under umbral composition. From our properties we derive new and interesting identities of Sheffer sequences of special polynomials for the powers of Sheffer pairs under umbral composition.
\end{abstract}

\section{Introduction and Preliminaries}
For $\alpha\in\mathbf{R}$, the Bernoulli polynomials of order $\alpha$ are defined by the generating function to be
\begin{equation}\label{eq:1}
\left(\frac{t}{e^{t}-1}\right)^{\alpha}e^{xt}=\sum_{n=0}^{\infty}B^{(\alpha)}_{n}(x)\frac{t^{n}}{n!},\quad \left(\text{see}\,\,\lbrack1,3,5\rbrack\right).
\end{equation}
In the special case, $x=0$, $B^{(\alpha)}_{n}(0)=B^{(\alpha)}_{n}$ are called the $n$-th Bernoulli numbers of order $\alpha$.\\
The Stirling number of the first kind is defined by
\begin{equation}\label{eq:2}
(x)_n=\sum_{k=0}^{n}S_{1}(n,k)x^{k},
\end{equation}
where $(x)_{n}=x(x-1)\cdots(x-n+1)$.\\
From (\ref{eq:2}), we note that
\begin{equation}\label{eq:3}
x^{(n)}=x(x+1)\cdots(x+n-1)=(-1)^{n}(-x)_{n}=\sum_{k=0}^{n}\left\vert S_{1}(n,k)\right\vert x^{k}.
\end{equation}
Let $\mathcal{F}$ be the set of all formal power series in the variable $t$ over $\mathbf{C}$ with
\begin{equation}\label{eq:4}
\mathcal{F}=\left\{ f(t)=\sum_{k=0}^{\infty}\frac{a_{k}}{k!}t^{k} \vert a_{k}\in\mathbf{C}\right\}.
\end{equation}
Suppose that $\mathbb{P}$ is the algebra of polynomials in the variable $x$ over $\mathbf{C}$ and $\mathbb{P}^{*}$ is the vector space of all linear functionals on $\mathbb{P}$. The action of the linear functional $L$ on a polynomial $p(x)$ is denoted by $\langle L\vert p(x)\rangle$. For $f(t)=\sum_{k=0}^{\infty}\frac{a_{k}}{k!}t^{k}\in\mathcal{F}$, let us define a linear functional on $\mathbb{P}$ by setting
\begin{equation}\label{eq:5}
\langle f(t)\vert x^{n}\rangle=a_{n},\quad (n\geq 0),\,\,\,\,(\text{see}\,\, \lbrack 2, 4\rbrack).
\end{equation}
By (\ref{eq:4}) and (\ref{eq:5}), we easily get
\begin{equation}\label{eq:6}
\langle t^{k}\vert x^{n} \rangle = n!\delta_{n,k}\quad (n,k\geq 0),\,\,\,\,(\text{see}\,\, \lbrack 2, 4\rbrack),
\end{equation}
where $\delta_{n,k}$ is the Kronecker's symbol.\\
For $f_{L}(t)=\sum_{k=0}^{\infty}\frac{\langle L\vert x^{k} \rangle}{k!}t^{k}$, we have $\langle f_{L}(t)\vert x^{n} \rangle=\langle L\vert x^{n} \rangle$.\\
Thus, we note that the map $L\longmapsto f_{L}(t)$ is a vector space isomorphism from $\mathbb{P}^{*}$ onto $\mathcal{F}$. Henceforth, $\mathcal{F}$ is thought of as both a formal power series and a linear functional. We call $\mathcal{F}$ the umbral algebra. The umbral calculus is the study of umbral algebra (see $\lbrack 4 \rbrack$).\\
The order $O(f(t))$ of the nonzero power series $f(t)$ is the smallest integer $k$ for which the coefficient of $t^{k}$ does not vanish (see $\lbrack2,4\rbrack$).\\
If $O(f(t))=0$, then $f(t)$ is called an invertible series. If $O(f(t))=1$, then $f(t)$ is called a delta series. For $O(f(t))=1$ and $O(g(t))=0$, there exists a unique sequence $s_{n}(x)$ of polynomials such that $\langle g(t)f(t)^{k}\vert s_{n}(x) \rangle = n!\delta_{n,k}$ for $n,k\geq 0$.\\
The sequence $s_{n}(x)$ is called the Sheffer sequence for $(g(t),f(t))$ which is denoted by $s_{n}(x)\sim(g(t),f(t))$.\\
Let $f(t)\in\mathcal{F}$ and $p(x)\in\mathbb{P}$. Then we see that
\begin{equation}\label{eq:7}
f(t)=\sum_{k=0}^{\infty}\frac{\langle f(t)\vert x^{k} \rangle}{k!}t^{k},\quad p(x)=\sum_{k=0}^{\infty}\frac{\langle t^k\vert p(x) \rangle}{k!}x^{k},\quad (\text{see}\,\, \lbrack 4\rbrack).
\end{equation}
By (\ref{eq:7}), we easily see that
\begin{equation}\label{eq:8}
t^{k}p(x)=p^{(k)}(x)=\frac{d^{k}p(x)}{dx^{k}},\quad (\text{see}\,\, \lbrack 2, 4\rbrack).
\end{equation}
Let $s_{n}(x)\sim(g(t),f(t))$. Then the generating function of Sheffer sequence $s_{n}(x)$ is given by
\begin{equation}\label{eq:9}
\frac{1}{g(\bar{f}(t))}e^{x\bar{f}(t)}=\sum_{k=0}^{\infty}s_{k}(x)\frac{t^{k}}{k!},\quad (\text{see}\,\, \lbrack 2, 4\rbrack),
\end{equation}
where $\bar{f}(t)$ is the compositional inverse of $f(t)$.\\
For $p_{n}(x)\sim(1,f(t))$, $q_{n}(x)\sim(1,g(t))$, we note that
\begin{equation}\label{eq:10}
q_{n}(x)=x\left(\frac{f(t)}{g(t)}\right)^{n}x^{-1}p_{n}(x),\quad (\text{see}\,\, \lbrack 2, 4\rbrack).
\end{equation}
The pair $(g(t),f(t))$ will be called a Sheffer pair where $O(g(t))=0$ and $O(f(t))=1$ (\text{see}\,\, \lbrack 2, 4\rbrack). Let $m$ be nonnegative integer. The $m$-th power of an invertible series is denoted by $(g(t))^{m}$, while the compositional power of a delta series $f(t)$ is denoted by $f^{m}(t)=\underbrace{f\circ f \circ \cdots \circ f}_{m-times}(t)$. Let $p_{n}(x)$ and $q_{n}(x)=\sum_{k=0}^{n}q_{n,k}x^{k}$ be sequences of polynomials. Then the umbral composition of $q_{n}(x)$ with $p_{n}(x)$ is defined by
\begin{equation}\label{eq:11}
\left(q_{n}\circ p\right)(x)=\sum_{k=0}^{n}q_{n,k}p_{k}(x),\quad (\text{see}\,\, \lbrack 2, 4\rbrack).
\end{equation}
Suppose that $s_{n}(x)\sim (g(t),f(t))$ and $r_{n}(x)\sim (h(t),l(t))$\\
Then we note that
\begin{equation}\label{eq:12}
(r_{n}\circ s)(x)=r_{n}(s(x))\sim (g(t)h(f(t)),l(f(t))).
\end{equation}
The identity under umbral composition is the sequence $x^{n}$ and the inverse of sequence $s_{n}(x)$ is the Sheffer sequence for $\left(g(\bar{f}(t))^{-1},\bar{f}(t)\right)$ (\text{see}\,\, \lbrack 2, 4\rbrack).\\
By (\ref{eq:12}), we easily see that the $m$-th power under umbral composition of $s_{n}(x)\sim (g(t),f(t))$ is given by
\begin{equation}\label{eq:13}
s^{(m)}_{n}(x)\sim\left(\prod_{i=0}^{m-1}g(f^{i}(t)),f^{m}(t)\right),\quad \text{where}\,\, m\in\mathbf{N}.
\end{equation}
For $n\geq 0$, let us assume that
\begin{equation}\label{eq:14}
s_{n}(x)=\sum_{k=0}^{n}s_{n,k}x^{k}=\sum_{k=0}^{\infty}s_{n,k}x^{k},
\end{equation}
where we agree that $s_{i,j}=0$ if $i<j$.\\
If we define $s^{(m)}_{n}(x)$ by
\begin{equation}\label{eq:15}
s^{(m)}_{n}(x)=\sum_{k=0}^{n}s^{(m)}_{n,k}x^{k}=\sum_{k=0}^{\infty}s^{(m)}_{n,k}x^{k},
\end{equation}
then, by (\ref{eq:11}),(\ref{eq:14}) and (\ref{eq:15}), we easily get
\begin{equation}\label{eq:16}
s^{(m)}_{n,k}=\sum_{l_{1},\cdots,l_{m-1}=0}^{n}s_{n,l_{1}}s_{l_{1},l_{2}}\cdots s_{l_{m-2},l_{m-1}}s_{l_{m-1},k},\quad (\text{see}\,\, \lbrack 2\rbrack).
\end{equation}
From (\ref{eq:9}) and (\ref{eq:13}), we can derive the generating function of $s^{(m)}_{n}(x)$ as follows:
\begin{align}\label{eq:17}
\sum_{k=0}^{\infty}\frac{s^{(m)}_{k}(x)}{k!}t^{k}&=\left(\frac{1}{\prod_{i=0}^{m-1}g(f^{i}(\bar{f}^{m}(t)))} \right)e^{x\bar{f}^{m}(t)}\\&=\left(\prod_{i=0}^{m-1}g(\bar{f}^{(m-i)}(t))\right)^{-1}e^{x\bar{f}^{m}(t)}\,.\nonumber
\end{align}
In this paper, we study some properties of Sheffer sequences for the powers of Sheffer pairs under umbral composition. From our properties, we derive new and interesting identities of Sheffer sequences of special polynomials for the powers of Sheffer pairs under umbral composition.

\section{Some identities of special polynomials.}
Let us take the sequence $s_{n}(x)$ of special polynomial as follows:
\begin{equation}\label{eq:18}
s_{n}(x)=x^{(n)}=\sum_{k=0}^{n}\vert S_{1}(n,k)\vert x^{k}\sim (1,f(t)=1-e^{-t}).
\end{equation}
For $m\in\mathbf{N}$, let us assume that the $m$-th power under umbral composition of $s_{n}(x)$ is given by
\begin{equation}\label{eq:19}
s^{(m)}_{n}(x)=\sum_{k=0}^{n}s^{(m)}_{n,k}x^{k}.
\end{equation}
By (\ref{eq:16}), (\ref{eq:18}) and (\ref{eq:19}), we get
\begin{align}\label{eq:20}
s^{(m)}_{n,k}&=\sum_{l_{1},\cdots,l_{m-1}=0}^{n}\vert S_{1}(n,l_{1})\vert\vert S_{1}(l_{1},l_{2})\vert\cdots\vert S_{1}(l_{m-1},k)\vert\\
&=\sum_{l_{1},\cdots,l_{m-1}=0}^{n}\vert S_{1}(n,l_{1})S_{1}(l_{1},l_{2})\cdots S_{1}(l_{m-1},k)\vert\,.\nonumber
\end{align}
It is known that
\begin{equation}\label{eq:21}
x^{n}\sim (1,t),\quad s_{n}(x)=x^{(n)}\sim (1,f(t)=1-e^{-t}).
\end{equation}
By (\ref{eq:10}) and (\ref{eq:21}), we get
\begin{equation}\label{eq:22}
s_{n}(x)=x\left(\frac{t}{f(t)}\right)^{n}x^{-1}x^{n}=x\left(\frac{t}{f(t)}\right)^{n}x^{n-1}.
\end{equation}
From (\ref{eq:22}), we note that
\begin{align}\label{eq:23}
f(t)^{m}x^{-1}s_{n}(x)&=f(t)^{m}\left(\frac{t}{f(t)}\right)^{n}x^{n-1}=\left(\frac{t}{f(t)}\right)^{n-m}t^{m}x^{n-1}\\
&=\left(\frac{t}{1-e^{-t}}\right)^{n-m}t^{m}x^{n-1}=\sum_{l=0}^{\infty}\frac{(-1)^{l}B^{(n-m)}_{l}}{l!}t^{l+m}x^{n-1}\nonumber\\
&=\sum_{l=0}^{n-1-m}\frac{(-1)^{l}B^{(n-m)}_{l}}{l!}\left(n-1\right)_{l+m}x^{n-1-l-m}\, ,\nonumber
\end{align}
where $n\geq 1$, $0\leq m\leq n-1$.\\
For $n\geq 1$, by (\ref{eq:13}), (\ref{eq:18}), we get
\begin{align}\label{eq:24}
s^{(2)}_{n}(x)&=x\left(\frac{f(t)}{f^{2}(t)}\right)^{n}x^{-1}s_{n}(x)=x\left(\frac{f(t)}{1-e^{-f(t)}}\right)^{n}x^{-1}s_{n}(x)\\
&=x\sum_{k_{2}=0}^{n-1}\frac{B^{(n)}_{k_{2}}}{k_{2}!}(-1)^{k_{2}}f(t)^{k_{2}}x^{-1}s_{n}(x).\nonumber
\end{align}
From (\ref{eq:23}) and (\ref{eq:24}), we can derive the following equation (\ref{eq:25}):
\begin{align}\label{eq:25}
s^{(2)}_{n}(x)&=x\sum_{k_{2}=0}^{n-1}\frac{B^{(n)}_{k_{2}}}{k_{2}!}(-1)^{k_{2}}\sum_{k_{1}=0}^{n-1-k_{2}}\frac{(-1)^{k_{1}}B^{(n-k_{2})}_{k_{1}}}{k_{1}!}(n-1)_{k_{1}+k_{2}}x^{n-1-k_{1}-k_{2}}\\
&=\sum_{k_{2}=0}^{n-1}\sum_{k_{1}=0}^{n-1-k_{2}}\frac{(n-1)!(-1)^{k_{1}+k_{2}}}{k_{1}!k_{2}!(n-k_{1}-k_{2}-1)!}B^{(n)}_{k_{2}}B^{(n-k_{2})}_{k_{1}}x^{n-k_{1}-k_{2}}\nonumber\\
&=\sum_{k_{1}+k_{2}+l=n-1}^{}\left(\begin{array}{c}n-1\\k_{1},k_{2},l\end{array} \right)(-1)^{k_{1}+k_{2}}B^{(n)}_{k_{2}}B^{(n-k_{2})}_{k_{1}}x^{l+1}\nonumber\\
&=\sum_{k=1}^{n}\left\{\sum_{k_{1}+k_{2}=n-k}^{}\left(\begin{array}{c}n-1\\k_{1},k_{2},k-1\end{array}\right)(-1)^{k_{1}+k_{2}}B^{(n)}_{k_{2}}B^{(n-k_{2})}_{k_{1}} \right\}x^{k}.\nonumber
\end{align}
From $s^{(3)}_{n}(x)\sim (1,f^{3}(t))$ and $s^{(2)}(x)\sim (1,f^{2}(t))$, we get
\begin{align}\label{eq:26}
s^{(3)}_{n}(x)&=x\left(\frac{f^{2}(t)}{f^{3}(t)}\right)^{n}x^{-1}s^{(2)}_{n}(x)=x\left(\frac{f^{2}(t)}{1-e^{-f^{2}(t)}}\right)^{n}x^{-1}s^{(2)}_{n}(x)\\
&=x\sum_{k_{3}=0}^{n-1}\frac{B^{(n)}_{k_{3}}}{k_{3}!}(-1)^{k_{3}}\left(\frac{f(t)}{f^{2}(t)}\right)^{n-k_{3}}\left(f(t)\right)^{k_{3}}x^{-1}s_{n}(x)\nonumber\\
&=x\sum_{k_{3}=0}^{n-1}\frac{B^{(n)}_{k_{3}}}{k_{3}!}(-1)^{k_{3}}\left(\frac{f(t)}{1-e^{-f(t)}}\right)^{n-k_{3}}\left(f(t)\right)^{k_{3}}x^{-1}s_{n}(x).\nonumber
\end{align}
From (\ref{eq:1}), (\ref{eq:23}) and (\ref{eq:26}), we have
\begin{align}\label{eq:27}
s^{(3)}_{n}(x)&=x\sum_{k_{3}=0}^{n-1}\sum_{k_{2}=0}^{n-1-k_{3}}(-1)^{k_{2}+k_{3}}\frac{B^{(n)}_{k_{3}}B^{(n-k_{3})}_{k_{2}}}{k_{3}!k_{2}!}\\
&\quad\times \sum_{k_{1}=0}^{n-1-k_{3}-k_{2}}\frac{(-1)^{k_{1}}B^{(n-k_{2}-k_{3})}_{k_{1}}}{k_{1}!}(n-1)_{k_{1}+k_{2}+k_{3}}x^{n-1-k_{1}-k_{2}-k_{3}}\nonumber\\
&=\sum_{k_{1}+k_{2}+k_{3}+l=n-1}^{}\left(-1\right)^{k_{1}+k_{2}+k_{3}}\left(\begin{array}{c}n-1\\k_{1},k_{2},k_{3},l\end{array} \right)B^{(n)}_{k_{3}}B^{(n-k_{3})}_{k_2}B^{(n-k_{3}-k_{2})}_{k_1}x^{l+1}\nonumber\\
&=\sum_{k=1}^{n}\Bigg\lbrace\sum_{k_{1}+k_{2}+k_{3}=n-k}^{}\left(-1\right)^{k_{1}+k_{2}+k_{3}}\left(\begin{array}{c}n-1\\k_{1},k_{2},k_{3},k-1\end{array}\right)\nonumber\\
&\quad\times B^{(n)}_{k_{3}}B^{(n-k_{3})}_{k_2}B^{(n-k_{3}-k_{2})}_{k_1} \Bigg\rbrace x^k.\nonumber
\end{align}
Continuing this process, we get
\begin{align}\label{eq:28}
s^{(m)}_{n}(x)&=\sum_{k=1}^{n}\Bigg\lbrace\sum_{k_{1}+\cdots+k_{m}=n-k}^{}(-1)^{k_{1}+\cdots+k_{m}}\left(\begin{array}{c}n-1\\k_{1},\cdots,k_{m},k-1\end{array} \right)B^{(n)}_{k_{m}}\\
&\quad\times B^{(n-k_{m})}_{k_{m-1}}\cdots B^{(n-k_{m}-\cdots-k_{2})}_{k_1} \Bigg\rbrace x^{k}.\nonumber
\end{align}
Therefore, by (\ref{eq:19}), (\ref{eq:20}) and (\ref{eq:28}), we obtain the following theorem.
\begin{thm}\label{theorem1}
For $m, n\geq 1$, we have
\begin{align*}
&\sum_{l_{1},\cdots,l_{m-1}=0}^{n}\left\vert S_{1}(n,l_{1})S_{1}(l_{1},l_{2})\cdots S_{1}(l_{m-1},k)\right\vert\\
&=\sum_{k_{1}+\cdots+k_{m}=n-k}^{}(-1)^{k_{1}+\cdots+k_{m}}\left(\begin{array}{c}n-1\\k_{1},\cdots,k_{m},k-1\end{array}\right)\\
&~~~~~~~~~~~~~~~~\quad\times B^{(n)}_{k_{m}}B^{(n-k_{m})}_{k_{m-1}}\cdots B^{(n-k_{m}-k_{m-1}-\cdots-k_{2})}_{k_1}.
\end{align*}
\end{thm}
Let us consider the following Sheffer sequence:
\begin{equation}\label{eq:29}
s_{n}(x)=L_{n}(x)=\sum_{k=0}^{n}L(n,k)(-x)^{k}\sim\left(1,f(t)=\frac{t}{t-1}\right),
\end{equation}
where $L(n,k)$ are the Lah numbers with
\begin{align}\label{eq:30}
&L(n,k)=\left(\begin{array}{c}n-1\\k-1\end{array}\right)\frac{n!}{k!},\quad \text{for}\,\, 1\leq k \leq n,\\
&L(n,k)=0,\quad \text{for}\,\, k> n \geq 1,\nonumber\\
&L(n,0)=0,\quad \text{for}\,\, n\geq 1,\nonumber\\
&L_{n}(0,0)=1.\nonumber
\end{align}
For $n\geq 1$, $0\leq m \leq n-1$, we have
\begin{align}\label{eq:31}
f(t)^{m}x^{-1}s_{n}(x)&=f(t)^{m}\left(\frac{t}{f(t)}\right)^{n}x^{n-1}=\left(\frac{t}{f(t)}\right)^{n-m}t^{m}x^{n-1}\\
&=(t-1)^{n-m}t^{m}x^{n-1}\nonumber\\
&=\sum_{l=0}^{n-m-1}\left(\begin{array}{c}n-m\\l\end{array}\right)(-1)^{n-m-l}(n-1)_{m}t^{l}x^{n-1-m}\nonumber\\
&=\sum_{l=0}^{n-m-1}\left(\begin{array}{c}n-m\\l\end{array}\right)(-1)^{n-m-l}(n-1)_{m}(n-1-m)_{l}x^{n-1-m-l}\nonumber\\
&=\sum_{l=0}^{n-m-1}\left(\begin{array}{c}n-m\\l\end{array}\right)(-1)^{n-m-l}(n-1)_{l+m}x^{n-1-m-l}.\nonumber
\end{align}
For $n \geq 1$, from $s^{(2)}_{n}(x)\sim\left(1,f^{2}(t)\right)$ and $s_{n}(x)\sim\left(1,f(t)=\frac{t}{t-1}\right)$, we get
\begin{align}\label{eq:32}
s^{(2)}_{n}(x)&=x\left(\frac{f(t)}{f^{2}(t)}\right)^{n}x^{-1}s_{n}(x)=x\left(f(t)-1\right)^{n}x^{-1}s_{n}(x)\\
&=x\sum_{k_{2}=0}^{n-1}\left(\begin{array}{c}n\\k_{2}\end{array}\right)(-1)^{n-k_{2}}f(t)^{k_{2}}x^{-1}s_{n}(x).\nonumber
\end{align}
From (\ref{eq:31}) and (\ref{eq:32}), we can derive the following equation:
\begin{align}\label{eq:33}
s^{(2)}_{n}(x)&=x\sum_{k_{2}=0}^{n-1}\left(\begin{array}{c}n\\k_{2}\end{array}\right)(-1)^{n-k_{2}}\sum_{k_{1}=0}^{n-1-k_{2}}\left(\begin{array}{c}n-k_{2}\\k_{1}\end{array}\right)(-1)^{n-k_{2}-k_{1}}\\
&\qquad \times(n-1)_{k_{1}+k_{2}}x^{n-1-k_{1}-k_{2}}\nonumber\\
&=\sum_{k_{2}=0}^{n-1}\sum_{k_{1}=0}^{n-1-k_{2}}(-1)^{n-k_{2}+(n-k_{1}-k_{2})}\frac{n!}{(n-k_{1}-k_{2})!}\nonumber\\
&\qquad\times\left(\begin{array}{c}n-1\\k_{1},k_{2},n-1-k_{1}-k_{2}\end{array}\right)x^{n-k_{1}-k_{2}}\nonumber\\
&=\sum_{k_{1}+k_{2}+l=n-1}^{}(-1)^{n-k_{2}+l+1}\frac{n!}{(l+1)!}\left(\begin{array}{c}n-1\\k_{1},k_{2},l\end{array}\right)x^{l+1}\nonumber\\
&=\sum_{k=1}^{n}\left\{\sum_{k_{1}+k_{2}=n-k}^{} (-1)^{(n-k_{2})+k}\frac{n!}{k!}\left(\begin{array}{c}n-1\\k_{1},k_{2},k-1\end{array}\right)\right\}x^{k}\nonumber
\end{align}
From $s^{(3)}_{n}(x)\sim\left(1,f^{3}(t)\right)$ and $s^{(2)}_{n}(x)\sim\left(1,f^{2}(t)\right)$, we get
\begin{align}\label{eq:34}
s^{(3)}_{n}(x)&=x\left(\frac{f^{2}(t)}{f^{3}(t)}\right)^{n}x^{-1}s^{(2)}_{n}(x)=x\left(\frac{f^{2}(t)}{\frac{f^{2}(t)}{f^{2}(t)-1}}\right)^{n}x^{-1}s^{(2)}_{n}(x)\\
&=x\left(f^{2}(t)-1\right)^{n}x^{-1}s^{(2)}_{n}(x)=x\sum_{k_{3}=0}^{n-1}\left(\begin{array}{c}n\\k_{3}\end{array}\right)(-1)^{n-k_{3}}\left(f^{2}(t)\right)^{k_{3}}x^{-1}s^{(2)}_{n}(x)\nonumber\\
&=x\sum_{k_{3}=0}^{n-1}\left(-1\right)^{n-k_{3}}\left(\begin{array}{c}n\\k_{3}\end{array}\right)\left(\frac{f(t)}{f^{2}(t)}\right)^{n-k_{3}}f(t)^{k_{3}}x^{-1}s_{n}(x)\nonumber\\
&=x\sum_{k_{3}=0}^{n-1}\left(\begin{array}{c}n\\k_{3}\end{array}\right)\left(-1\right)^{n-k_{3}}\sum_{k_{2}=0}^{n-1-k_{3}}\left(\begin{array}{c}n-k_{3}\\k_{2}\end{array}\right)\left(-1\right)^{n-k_{2}-k_{3}}\left(f(t)\right)^{k_{2}+k_{3}}x^{-1}s_{n}(x).\nonumber
\end{align}
From (\ref{eq:31}) and (\ref{eq:34}), we have
\begin{align}\label{eq:35}
s^{(3)}_{n}(x)&=\sum_{k_{1}+k_{2}+k_{3}+l=n-1}^{}\left(-1\right)^{(n-k_{3})+(n-k_{2}-k_{3})+(l+1)}\frac{n!}{(l+1)!}\left(\begin{array}{c}n-1\\k_{1},k_{2},k_{3},l\end{array}\right) x^{l+1}\\
&=\sum_{k=1}^{n}\left\{\sum_{k_{1}+k_{2}+k_{3}=n-k}^{}\left(-1\right)^{(n-k_{3})+(n-k_{2}-k_{3})+k}\frac{n!}{k!}\left(\begin{array}{c}n-1\\k_{1},k_{2},k_{3},k-1\end{array}\right) \right\}x^k.\nonumber
\end{align}
Continuing this process, we get
\begin{align}\label{eq:36}
s^{(m)}_{n}(x)&=\sum_{k=1}^{n}\Bigg\lbrace\sum_{k_{1}+\cdots+k_{m}=n-k}^{}\left(-1\right)^{(n-k_{m})+\cdots+(n-k_{m}-k_{m-1}-\cdots-k_{2})+k}\frac{n!}{k!}\nonumber\\
&\quad\times\left(\begin{array}{c}n-1\\k_{1},k_{2},\cdots,k_{m},k-1\end{array}\right) \Bigg\rbrace x^{k}\\
&=\sum_{k=1}^{n}s^{(m)}_{n,k}x^{k},\quad \text{where}\,\, m\geq 1.\nonumber
\end{align}
By (\ref{eq:14}), (\ref{eq:15}), (\ref{eq:16}), (\ref{eq:29}) and (\ref{eq:36}), we easily get
\begin{align}\label{eq:37}
s^{(m)}_{n,k}&=\sum_{l_{1},\cdots,l_{m-1}=0}^{n}s_{n,l_{1}}s_{l_{1},l_{2}}\cdots s_{l_{m-1},k}\\
&=\sum_{l_{1},\cdots,l_{m-1}=0}^{n}\left(-1\right)^{l_{1}+l_{2}+\cdots+l_{m-1}+k}L(n,l_{1})L(l_{1},l_{2})\cdots L(l_{m-1},k).\nonumber
\end{align}
Therefore, by (\ref{eq:36}) and (\ref{eq:37}), we obtain the following theorem.
\begin{thm}\label{theorem2}
For $m,n\geq 1$, $1\leq k\leq n$, we have
\begin{align*}
&\sum_{l_{1},l_{2}\cdots,l_{m-1}=0}^{n}\left(-1\right)^{l_{1}+l_{2}+\cdots+l_{m-1}+k}L(n,l_{1})L(l_{1},l_{2})\cdots L(l_{m-1},k)\\
&=\sum_{k_{1}+\cdots,k_{m}=n-k}^{}\left(-1\right)^{(n-k_{m})+(n-k_{m}-k_{m-1})+\cdots+(n-k_{m}-\cdots-k_{2})+k}\frac{n!}{k!}\left(\begin{array}{c}n-1\\k_{1},k_{2},\cdots,k_{m},k-1\end{array}\right)
\end{align*}
\end{thm}

Let us take Abel sequence as follows:
\begin{align}\label{eq:38}
s_{n}(x)&=A_{n}(x:a)=x(x-an)^{n-1}=\sum_{k=1}^{n}\left(\begin{array}{c}n-1\\k-1\end{array}\right)(-an)^{n-k}x^{k}\\
&\sim\left(1,f(t)=te^{at}\right),\quad \text{where}\,\, a\not=0.\nonumber
\end{align}
Thus by (\ref{eq:38}), we get
\begin{equation}\label{eq:39}
s_{n,k}=\left(\begin{array}{c}n-1\\k-1\end{array}\right)\left(-an\right)^{n-k},\quad (n, k\geq 0)
\end{equation}
From (\ref{eq:16}) and (\ref{eq:39}), we note that
\begin{align}\label{eq:40}
s^{(m)}_{n,k}&=\sum_{l_{1},\cdots,l_{m-1}=0}^{n}s_{n,l_{1}}s_{l_{1},l_{2}}\cdots s_{l_{m-2},l_{m-1}}s_{l_{m-1},k}\\
&=\sum_{l_{1},\cdots,l_{m-1}=0}^{n}\left(\begin{array}{c}n-1\\l_{1}-1\end{array}\right)\left(\begin{array}{c}l_{1}-1\\l_{2}-1\end{array}\right)\cdots\left(\begin{array}{c}l_{m-2}-1\\l_{m-1}-1\end{array}\right)\left(\begin{array}{c}l_{m-1}-1\\k-1\end{array}\right)\nonumber\\
&\quad\times(-a)^{n-k}n^{n-l_{1}}l_{1}^{l_{1}-l_{2}}\cdots l_{m-2}^{l_{m-2}-l_{m-1}}l_{m-1}^{l_{m-1}-k}.\nonumber
\end{align}
From $s_{n}(x)=A_{n}\left(x:a\right)\sim\left(1,f(t)=te^{at}\right)$ and $x^{n}\sim \left(1,t\right)$, we note that
\begin{align}\label{eq:41}
f(t)^{m}x^{-1}s_{n}(x)&=f(t)^{m}\left(\frac{t}{f(t)}\right)^{n}x^{n-1}=\left(\frac{t}{f(t)}\right)^{n-m}t^{m}x^{n-1}\\
&=\left(\frac{t}{te^{at}}\right)^{n-m}t^{m}x^{n-1}=e^{-a(n-m)t}t^{m}x^{n-1}\nonumber\\
&=\sum_{l=0}^{n-1-m}\left(-a(n-m)\right)^{l}\frac{(n-1)_{l+m}}{l!}x^{n-1-l-m}.\nonumber
\end{align}
For $n\geq 1$, from $s^{(2)}_{n}(x)\sim\left(1,f^{2}(t)\right)$ and $s_{n}(x)\sim\left(1,f(t)=te^{at}\right)$, we get
\begin{align}\label{eq:42}
s^{(2)}_{n}(x)&=x\left(\frac{f(t)}{f^{2}(t)}\right)^{n}x^{-1}s_{n}(x)=x\left(\frac{f(t)}{f(t)e^{af(t)}}\right)^{n}x^{-1}s_{n}(x)\\
&=xe^{-anf(t)}x^{-1}s_{n}(x)=x\sum_{k_{2}=0}^{n-1}\frac{(-an)^{k_2}}{k_{2}!}\left(f(t)\right)^{k_2}x^{-1}s_{n}(x).\nonumber
\end{align}
From (\ref{eq:41}) and (\ref{eq:42}), we can derive the following equation (\ref{eq:43}):
\begin{align}\label{eq:43}
s^{(2)}_{n}(x)&=\sum_{k_{2}=0}^{n-1}\sum_{k_{1}=0}^{n-1-k_{2}}\left(\begin{array}{c}n-1\\k_{1},k_{2},n-1-k_{1}-k_{2}\end{array}\right)\left(-an\right)^{k_{2}}\left(-a(n-k_{2})\right)^{k_{1}}x^{n-k_{1}-k_{2}}\\
&=\sum_{k_{1}+k_{2}+l=n-1}\left(\begin{array}{c}n-1\\k_{1},k_{2},l\end{array}\right)\left(-an\right)^{k_{2}}\left(-a(n-k_{2})\right)^{k_{1}}x^{l+1}\nonumber\\
&=\sum_{k=1}^{n}\left\{\sum_{k_{1}+k_{2}=n-k}^{}\left(\begin{array}{c}n-1\\k_{1},k_{2},k-1\end{array}\right)\left(-an\right)^{k_{2}}\left(-a(n-k_{2})\right)^{k_{1}} \right\}x^{k}.\nonumber
\end{align}
From $s^{(3)}_{n}(x)\sim\left(1,f^{3}(t)\right)$ and $s^{(2)}_{n}(x)\sim\left(1,f^{2}(t)\right)$, we get
\begin{align}\label{eq:44}
s^{(3)}_{n}(x)&=x\left(\frac{f^{2}(t)}{f^{3}(t)}\right)x^{-1}s^{(2)}_{n}(x)=xe^{-anf^{2}(t)}x^{-1}s^{(2)}_{n}(x)\\
&=x\sum_{k_{3}=0}^{n-1}\frac{(-an)^{k_{3}}}{k_{3}!}\left(f^{2}(t)\right)^{k_{3}}x^{-1}s^{(2)}_{n}(x)\nonumber\\
&=x\sum_{k_{3}=0}^{n-1}\frac{(-an)^{k_{3}}}{k_{3}!}\left(\frac{f(t)}{f^{2}(t)}\right)^{n-k_{3}}f(t)^{k_{3}}x^{-1}s_{n}(x)\nonumber\\
&=x\sum_{k_{3}=0}^{n-1}\frac{(-an)^{k_{3}}}{k_{3}!}e^{-a(n-k_{3})f(t)}\left(f(t)\right)^{k_{3}}x^{-1}s_{n}(x)\nonumber\\
&=x\sum_{k_{3}=0}^{n-1}\frac{(-an)^{k_{3}}}{k_{3}!}\sum_{k_{2}=0}^{n-1-k_{3}}\frac{(-a(n-k_{3}))^{k_{2}}}{k_{2}!}\left(f(t)\right)^{k_{2}+k_{3}}x^{-1}s_{n}(x).\nonumber
\end{align}
From (\ref{eq:41}) and (\ref{eq:44}), we can derive the following equation (\ref{eq:45}):
\begin{align}\label{eq:45}
s^{(3)}_{n}(x)&=\sum_{k_{1}+k_{2}+k_{3}+l=n-1}\left(\begin{array}{c}n-1\\k_{1},k_{2},k_{3},l\end{array}\right)\left(-an\right)^{k_{3}}\left(-a(n-k_{3})\right)^{k_{2}}\left(-a(n-k_{2}-k_{3})\right)^{k_{1}}x^{l+1}\\
&=\sum_{k=1}^{n}\Bigg\lbrace \sum_{k_{1}+k_{2}+k_{3}=n-k}\left(\begin{array}{c}n-1\\k_{1},k_{2},k_{3},k-1\end{array}\right)\left(-an\right)^{k_{3}}\left(-a\left(n-k_{3}\right)\right)^{k_{2}}\nonumber\\
&\quad\times\left(-a\left(n-k_{2}-k_{3}\right)\right)^{k_{1}} \Bigg\rbrace x^{k}.\nonumber
\end{align}
Continuing this process, we get
\begin{align}\label{eq:46}
s^{(m)}_{n}(x)&=\sum_{k=1}^{n}\Bigg\lbrace\sum_{k_{1}+\cdots+k_{m}=n-k}\left(\begin{array}{c}n-1\\k_{1},k_{2},\cdots,k_{m},k-1\end{array}\right)\nonumber\\
&\quad\times\left(\prod_{i=1}^{m}\left(-a\left(n-k_{m}-\cdots-k_{i+1}\right)\right)^{k_{i}} \right)\Bigg\rbrace x^{k}.
\end{align}
Therefore, by (\ref{eq:40}) and (\ref{eq:46}), we obtain the following theorem.
\begin{thm}\label{theorem3}
For $n,m\geq 1$, $1\leq k\leq n$, we have
\begin{align*}
&\sum_{l_{1},\cdots,l_{m-1}=0}^{n}\left(\begin{array}{c}n-1\\l_{1}-1\end{array}\right)\left(\begin{array}{c}l_{1}-1\\l_{2}-1\end{array}\right)\cdots\left(\begin{array}{c}l_{m-2}-1\\l_{m-1}-1\end{array}\right)\left(\begin{array}{c}l_{m-1}-1\\k-1\end{array}\right)(-an)^{n-l_{1}}\\
&\quad\times(-al_{1})^{l_{1}-l_{2}}\cdots(-al_{m-2})^{l_{m-2}-l_{m-1}}(-al_{m-1})^{l_{m-1}-k}\\
&=\sum_{k_{1}+k_{2}+\cdots+k_{m}=n-k}^{}\left(\begin{array}{c}n-1\\k_{1},k_{2},\cdots,k_{m},k-1\end{array}\right)\left(\prod_{i=1}^{m}\left(-a\left(n-k_{m}-\cdots-k_{i+1}\right)\right)^{k_{i}} \right).
\end{align*}
\end{thm}
\noindent$\large{\mathbf{Remark}}$. Let us consider the Mittag-Leffler sequences as follows:
\begin{align}\label{eq:47}
s_{n}(x)=M_{n}(x)&=\sum_{r=0}^{n}\left(\begin{array}{c}n\\r\end{array}\right)\frac{(n-1)!}{(r-1)!}2^{r}(x)_{r}\\
&=\sum_{k=0}^{n}\left\{\sum_{r=k}^{n}\left(\begin{array}{c}n\\r\end{array}\right)\frac{(n-1)!}{(r-1)!}2^{r}S_{1}(r,k)\right\}x^{k}\nonumber\\
&\quad\sim\left(1,\frac{e^{t}-1}{e^{t}+1}=f(t)\right).\nonumber
\end{align}
By the same method, we get, for $m, n\geq 1$, $1\leq k\leq n$,
\begin{align*}
&\sum_{l_{1},\cdots,l_{m-1}=0}^{n}\sum_{r_{1}=l_{1}}^{n}\cdots\sum_{r_{m-1}=l_{m-1}}^{l_{m-2}}\sum_{r_{m}=k}^{l_{m-1}}
\left(\begin{array}{c}n\\r_{1}\end{array}\right)\left(\begin{array}{c}l_{1}\\r_{2}\end{array}\right)\cdots\left(\begin{array}{c}l_{m-2}\\r_{m-1}\end{array}\right)\left(\begin{array}{c}l_{m-1}\\r_{m}\end{array}\right)\\
&\quad\times\frac{(n-1)!(l_{1}-1)!\cdots (l_{m-2}-1)!(l_{m-1}-1)!}{(r_{1}-1)!(r_{2}-1)!\cdots(r_{m-1}-1)!(r_{m}-1)!}\times 2^{r_{1}+r_{2}+\cdots+r_{m}}\\
&\quad\times S_{1}(r_{1},l_{1})S_{1}(r_{2},l_{2})\cdots S_{1}(r_{m-1},l_{m-1})S_{1}(r_{m},k)\\
&=\sum_{k_{1}+\cdots+k_{2m}=n-k}^{}\left(\begin{array}{c}n-1\\k_{1},\cdots ,k_{2m},k-1\end{array}\right)\left(\prod_{i=0}^{m-1}E_{2i+1}^{(k_{1}+\cdots+k_{2i}-n)}B_{2i+2}^{(n-k_{1}-\cdots-k_{2i})} \right)\\
&\quad\times\left(\prod_{i=0}^{m-1}2^{n-(k_{1}+k_{2}+\cdots+k_{2i})} \right).
\end{align*}
Here, for $\alpha\in\mathbf{R}$, the Euler polynomials of order $\alpha$ are defined by the generating function to be
\begin{equation}\label{eq:1}
\left(\frac{2}{e^{t}+1}\right)^{\alpha}e^{xt}=\sum_{n=0}^{\infty}E^{(\alpha)}_{n}(x)\frac{t^{n}}{n!},\quad \left(\text{see}\,\,\lbrack1,3,5\rbrack\right).
\end{equation}
In the special case, $x=0$, $E^{(\alpha)}_{n}(0)=E^{(\alpha)}_{n}$ are called the $n$-th Euler numbers of order $\alpha$.\\

\author{Department of Mathematics, Sogang University, Seoul 121-742, Republic of Korea
\\e-mail: dskim@sogang.ac.kr}\\
\\
\author{Department of Mathematics, Kwangwoon University, Seoul 139-701, Republic of Korea
\\e-mail: tkkim@kw.ac.kr}

\begin{thebibliography}{99}
\bibitem{Acikgoz}
M. Acikgoz, D. Erdal, S. Araci,
\emph {A new approach to q-Bernoulli numbers and q-Bernoulli polynomials},
Adv. Difference Equ. 2010(2010), Art. ID 951764, 9pp.

\bibitem{KKLR}
D.S. Kim, T. Kim, S.H. Lee, C.S. Ryoo,
\emph{Sheffer sequences for the powers of Sheffer pairs under umbral composition},
Adv. Stud. Contemp. Math. 23 (2013), no. 2.

\bibitem{Kim}
T. Kim,
\emph{Symmetry p-adic invariant integral on $\mathbf{Z}_{p}$ for Bernoulli and Euler polynomials},
J. Difference Equ. Appl. 14(2008), no. 12, 1267--1277.

\bibitem{Roman}
S. Roman,
\emph{The umbral Calculus},
Dover Publ. Inc., New York, 1984.

\bibitem{Sen}
E. Sen,
\emph{Theorems on Apostol--Euler polynomials of higher order arising from Euler basis},
Adv. Stud. Contemp. Math. 23 (2013), no. 2.
\end{thebibliography}
\end{document}